\documentclass[a4paper]{amsart}

\usepackage[utf8]{inputenc}
\usepackage[margin=3.5cm]{geometry}
\usepackage{mathtools}
\usepackage{amssymb}
\usepackage[sortcites=true,giveninits=true,url=false,isbn=false]{biblatex}
\usepackage{hyperref}
\usepackage[capitalise]{cleveref}

\newtheorem*{theorem}{Theorem}
\newtheorem*{lemma}{Lemma}

\theoremstyle{definition}
\newtheorem*{definition}{Definition}

\newcommand{\diam}{\operatorname{diam}}
\newcommand{\card}{\operatorname{card}}

\renewcommand{\H}{\mathcal H}

\begin{document}
	
	\author{David Bate}
	\email{david.bate@warwick.ac.uk}
	\address{Zeeman Building, University of Warwick, Coventry CV4 7AL. ORCiD: \href{https://orcid.org/0000-0003-0808-2453}{0000-0003-0808-2453}}
	
	\begin{abstract}
		We prove that any Lipschitz map that satisfies a condition inspired by the work of David may be decomposed into countably many bi-Lipschitz pieces.
	\end{abstract}

	\title{Qualitative Lipschitz to bi-Lipschitz decomposition}
	
	
	\thanks{This work was supported by the European Union's Horizon 2020 research and innovation programme (Grant agreement No. 948021)}
	
	\maketitle

	Fix complete metric spaces $(X,d)$ and $(Y,\rho)$ and $s>0$.
	A map $f\colon X\to Y$ is \emph{Lipschitz} if there exists $L\geq 0$ such that $\rho(f(x),f(y))\leq L d(x,y)$ for all $x,y\in X$.
	An injective Lipschitz map is \emph{bi-Lipschitz} if $f^{-1}$ is Lipschitz.
	Let $\H^s$ denote the $s$-dimensional Hausdorff (outer) measure on a metric space and, for $S$ a subset of a metric space, let
	\begin{equation*}
			\Theta^{*,s}(S,x):=\limsup_{r\to 0}
		\frac{
			\H^s((B(x,r) \cap S))
		}{
			(2r)^s
		}
		\quad \text{and} \quad
			\Theta^{s}_*(S,x):=\liminf_{r\to 0}
		\frac{
			\H^s((B(x,r) \cap S))
		}{
			(2r)^s
		}
		.
	\end{equation*}

	\begin{theorem}
		Let $f\colon X\to Y$ be Lipschitz, $\H^s(X)<\infty$ and, for $\H^s$-a.e.\ $x\in X$, suppose that
		\begin{equation}\label{main-david-hyp}
			\limsup_{r\to 0}
			\frac{
				\H^s(B(f(x), \lambda_x r) \setminus f(B(x,r)))
			}{
				(2\lambda_x r)^s
			}
			< \frac{1}{2}\Theta_*^s(Y,f(x))
		\end{equation}
		for some $0<\lambda_x\leq 1$.
		Then there exists a countable Borel decomposition $X=N\cup \bigcup_i X_i$ with $\H^s(N)=0$ such that each $f|_{X_i}$ is bi-Lipschitz.
	\end{theorem}
	This theorem appeared in Bate and Li \cite[Theorem 1.2, iii) $\Rightarrow$ (R)]{bateli} as an intermediate step in order to prove characterisations of rectifiable subsets of a metric space.
	These characterisations initiated a line of research of geometric measure theory in metric spaces and form a cornerstone of the results in \cite{tan,proj}.
	The Theorem is a qualitative analogue of results of David \cite{david} and Semmes \cite{semmes}
	and the proof in \cite{bateli} adapts the 85 page argument in \cite{semmes} to the qualitative setting.
	We give a new self contained proof.

	A form of the following condition is present in \cite{semmes,david,bateli}.
	\begin{definition}
		For $0<\kappa < \lambda\leq 1$ and $0<\xi<1$, a function $f\colon V\subset X \to Y$ satisfies the condition $D(\lambda, \kappa,\xi)$ on a set $S\subset V$ if, for all $r < \diam S$,
		\begin{equation*}
			\H^s(B(f(x),\lambda r)\setminus f(V\cap B(x,r))) < \frac{1}{2} \xi\H^s(B(f(x),\kappa r)).
		\end{equation*}
	\end{definition}

	A key idea of David \cite{david} uses $D(\lambda,\kappa,\xi)$ to deduce bi-Lipschitz bounds on $f$.

	\begin{lemma}
		Let $0<\kappa < \lambda \leq 1$, $0<\xi<1$ and suppose $f\colon V\subset X \to Y$ satisfies $D(\lambda,\kappa,\xi)$ on $S\subset V$.
		Let $x,y\in S$, set $r=d(x,y)/4$ and suppose $\H^s(B(f(x),\kappa r))\geq \H^s(B(f(y),\kappa r))$.
		If
		\begin{equation}
			\label{david-assum-meas}
			\H^s(f(V\cap B(x,r)) \cap f(V\cap B(y,r)))
			\leq
			(1-\xi) \H^s\left(B\left(f(x),\kappa r\right)\right)
		\end{equation}
		then
		\begin{equation}
			\label{bilip-conclusion}
			\rho(f(x),f(y))\geq (\lambda-\kappa) \frac{d(x,y)}{4}.
		\end{equation}
	\end{lemma}

	\noindent\emph{Proof.}
		Suppose that \eqref{bilip-conclusion} does not hold.
		Then by the triangle inequality,
		\begin{equation*}
			B\left(f(x),\lambda r \right)\cap B\left(f(y), \lambda r \right) \supset B\left(f(x), \kappa r \right).
		\end{equation*}
		Combining this with $D(\lambda,\kappa,\xi)$ negates \eqref{david-assum-meas}.
		Indeed, it gives
		\begin{align*}
			\H^s(f(V\cap B(x,r))\cap f(V\cap B(y,r)))
			&>
			\H^s\left(B\left(f(x),\lambda r \right)\cap B\left(f(y), \lambda r \right)\right)
		      \\&\quad - \xi \H^s\left(B\left(f(x),\kappa r \right)\right)
			\\&\geq
			(1-\xi) \H^s\left(B\left(f(x),\kappa r\right)\right).
		\end{align*}

	\begin{proof}[Proof of the Theorem]

		First let $V\subset X$ be Borel and note that \eqref{main-david-hyp} holds for the function $f|_V$ and for $\H^s$-a.e.\ $x\in V$.
		Indeed, for $\H^s$-a.e.\ $x\in V$, $\Theta^{*,s}(X\setminus V,x)=0$ (see \cite[§2.10.18]{federer}) and for such an $x$,
		\begin{align}
			\label{D-subset}
			\limsup_{r\to 0}
			\frac{
				\H^s(B(f(x), \lambda_x r) \setminus f(V \cap B(x,r)))
			}{
				(2\lambda_x r)^s
			}
			\notag
			&\leq
			\limsup_{r\to 0}
			\frac{
				\H^s(B(f(x), \lambda_x r) \setminus f(B(x,r)))
			}{
				(2\lambda_x r)^s
			}
			\notag
			\\&\quad +
			\limsup_{r\to 0}
			\frac{
				L^s \H^s(B(x, r) \setminus V)
			}{
				(2\lambda_x r)^s
			}
			\notag
			\\&<
			\frac{1}{2}\Theta_*^s(Y,f(x)) +
			\frac{
			L^s
			}{
			\lambda_x^s}
			\Theta^{*,s}(X\setminus V,x).
		\end{align}
		The first inequality uses $\H^s(f(A))\leq L^s\H^s(A)$ for any $A\subset X$, see \cite[Corollary 2.10.11]{federer}.
		A consequence of \eqref{D-subset} is
		\begin{equation}
			\label{pos-meas}
			\H^s(V)>0 \Rightarrow \H^s(f(V))>0.
		\end{equation}

		Now let $V_1\subset X$ be Borel with $\H^s(V_1)>0$.
		By the coarea formula (see \cite[§2.10.25]{federer}), for $\H^s$-a.e.\ $y\in Y$, $\card f^{-1}(y)<\infty$.
		Hence, by \eqref{pos-meas}, for $\H^s$-a.e. $x\in V_1$,
		\begin{equation}
			\label{finite-card}
			\card\{x'\in V_1 : f(x')=f(x)\} <\infty.
		\end{equation}
		Let $V_2\subset V_1$ be a positive measure Borel set for which \eqref{finite-card} holds for all $x\in V_2$.
		By the Lusin-Novikov theorem (see \cite[Exercise 18.14]{kechris}
		\footnote{In fact, since $f$ is continuous, \eqref{finite-card} holds and $V_2$ may be chosen compact, the result is elementary.}
		), there exists a Borel function
			$g \colon f(V_2) \to V_2$
		such that $V_3:=g(f(V_2))$ is Borel and $f(g(y))=y$ for all $y\in f(V_2)$.
		By \eqref{pos-meas}, $\H^s(f(V_3)) = \H^s(f(V_2))>0$
		and hence, since $f$ is Lipschitz, $\H^s(V_3)>0$.

		Note that if $x\in X$ satisfies \eqref{main-david-hyp}, then it also satisfies \eqref{main-david-hyp} for all $0<\lambda\leq \lambda_x$.
		For $i\in \mathbb{N}$ let $1\geq \lambda_i \searrow 0$ and define $S_i$ to be the set of $x\in V_3$ for which
		\begin{equation*}
			\sup_{0<r<\lambda_i}
			\frac{
				\H^s(B(f(x), \lambda_i r) \setminus f(V_3 \cap B(x,r)))
			}{
				(2\lambda_i r)^s
			}
			< \inf_{0<r'<\lambda_i}
			\frac{1}{2} \left(1 - \lambda_i \right)^s
			\frac{
			\H^s(B(f(x),r'))
			}{(2r')^s}.
		\end{equation*}
		Then, by \eqref{D-subset}, the $S_i$ monotonically increase to a full measure subset of $V_3$.
		Therefore, there exist $i\in \mathbb{N}$ and $S'\subset S_i$ with $\H^s(S')>0$ and $\diam S' \leq \lambda_i$.
		For any $0<r<\lambda_i$, setting $r'=(1-\lambda_i^2)\lambda_i r$ shows that $f|_{V_3}$ satisfies $D(\lambda_i,(1-\lambda_i^2)\lambda_i,(1+\lambda_i)^{-s})$ on $S'$.
		Since $f|_{V_3}$ is injective, \eqref{david-assum-meas} holds for all $x,y\in S'$ and hence the Lemma implies that $f|_{S'}$ is bi-Lipschitz.

		The bi-Lipschitz condition extends to the closure of $S'$.
		Hence $\overline{S'} \cap V_1$ is a Borel subset of $V_1$ of positive measure on which $f$ is bi-Lipschitz.
		Since $V_1$ is an arbitrary Borel subset of positive measure, the conclusion follows by exhaustion (e.g.\ as in \cite[§3.2.14]{federer}).
	\end{proof}

	\printbibliography
\end{document}